\documentclass[reqno]{amsart}

\usepackage{amssymb,amsmath,amsthm,amsfonts,xcolor} 

\usepackage{comment,graphicx}

\usepackage{enumitem}
\usepackage{footnote}

\usepackage{color,setspace,wrapfig,multicol}

\usepackage[labelformat=empty]{subcaption}

\usepackage[mathscr]{eucal} 

\usepackage[colorlinks=true, pdfstartview=FitV, linkcolor=blue, citecolor=green]{hyperref}

\usepackage[left=3cm,right=3cm,top=2.5cm,bottom=2.5cm]{geometry}

\allowdisplaybreaks
\numberwithin{equation}{section}

\newtheorem{theorem}{Theorem}[section]

\theoremstyle{definition}

\newcommand{\be}{\mathbf{e}}

\newcommand{\R}{\mathbb{R}}
\newcommand{\cH}{\mathcal{H}}
\newcommand{\sC}{\mathscr{C}}
\renewcommand{\div}{\mathrm{div}}

\newcommand{\p}{\partial}
\newcommand{\inner}[2]{\left\langle#1,#2\right\rangle}

\setlist[itemize]{leftmargin=6mm} 

\begin{document}


\title[Minimality of Winterbottom shape]{On the minimality of the Winterbottom shape}

\author[Sh. Kholmatov] {Shokhrukh Yu. Kholmatov} 
\address[Sh. Kholmatov]{University of Vienna,  Oskar-Morgenstern Platz 1, 1090 Vienna 
(Austria)}
\email{shokhrukh.kholmatov@univie.ac.at}



\date{\today}

\begin{abstract}
In this short note we prove that the Winterbottom shape \cite{Winterbottom:1967} is a volume-constraint minimizer of the corresponding anisotropic capillary functional.
\end{abstract}

\maketitle

\section{Introduction} 

In this paper we study volume-constraint minimizers of the anisotropic capillary  functional in the upper half-space $\Omega:=\{x=(x_1,\ldots,x_n):\,\, x_n>0\}\subset \R^n:$
$$
\sC_{\Phi,\beta}(E):=P_\Phi(E,\Omega) - \beta \int_{\p\Omega}\chi_Ed\cH^{n-1},\quad E\in BV(\Omega;\{0,1\}),
$$
where $\Phi$ is an anisotropy -- a positively-one homogeneous convex function in $\R^n$ satisfying 
\begin{equation}\label{norm_bounds}
c_\Phi|x|\le \Phi(x)\le C_\Phi|x|,\quad x\in\R^n,
\end{equation}
for some $C_\Phi\ge c_\Phi>0$, 
$$
P_\Phi(E,\Omega) = \int_{\Omega\cap\p^*E} \Phi(\nu_E)d\cH^{n-1}
$$
is the $\Phi$-perimeter of $E$ in $\Omega,$ $\p^*E$ is the reduced boundary and $\nu_E$ is the generalized outer unit normal of $E,$ $\beta$ is constant -- a reative adhesion constant of $\p\Omega,$ and $\chi_E$ is the interior trace of $E$ along $\p\Omega,$ i.e.,
$$
\int_{\p\Omega}\chi_Ed\cH^{n-1} = \cH^{n-1}(\p\Omega\cap \p^*E).
$$
Recall that $E\in BV(\Omega;\{0,1\}) $ if and only if $E\in BV(\R^n;\{0,1\}),$ and in particular, $\chi_E\in L^1(\p\Omega)$ (see e.g. \cite{BKh:2018}). It is well-known that if $\beta\le -\Phi(-\be_n),$ where $\be_n:=(0,\ldots,0,1),$ then up to a translation the unique volume-constrained minimizer of $\sC_{\Phi,\beta}$ is a translation of the Wulff shape $W^\Phi:=\{\Phi^o\le 1\}$ of $\Phi$ in $\Omega,$ where 
$$
\Phi^o(x):=\max_{\Phi(y)=1}\,\,\inner{x}{y}
$$
is the dual anisotropy, where $\inner{\cdot}{\cdot}$ is the Euclidean scalar product. On the other hand, if $\beta\ge \Phi(\be_n),$ then one can readily check that 
$$
\inf_{E\in BV(\Omega;\{0,1\}),\,|E|=1}\,\,\sC_{\Phi,\beta}(E) = 
\begin{cases}
-\infty & \text{if $\beta>\Phi(\be_n),$}\\[1mm]
0 & \text{if $\beta=\Phi(\be_n),$}
\end{cases}
$$
and the minimum problem does not admit a solution. In case $\beta=0$ from  \cite[Theorem 1.3]{CRS:2016} we deduce the following relative isoperimetric inequality in $\Omega:$ 
\begin{equation}\label{crs_2016}
\frac{P_\Phi(E,\Omega)}{|E|^{\frac{n-1}{n}}} \ge  \frac{P_\Phi(W^\Phi,\Omega)}{|\Omega\cap W^\Phi|^{\frac{n-1}{n}}},\quad E\in BV(\Omega,\{0,1\}),
\end{equation}
for $0<|E\cap\Omega|<+\infty.$ It turns out (see \cite[p. 2979]{CRS:2016}) that the equality in \eqref{crs_2016} holds if and only if $E = b + rW^\Phi$ for some $b\in \p\Omega$ and $r>0,$ i.e., $E$ is a horizontal translation of scaled Wulff shapes. In particular, the set 
$
W_0^\Phi:=\Omega\cap W^\Phi
$
is a unique (up to a horizontal translation) solution to the minimum  problem 
$$
\inf_{E\in BV(\Omega;\{0,1\}),|E|=|W_0^\Phi|}\,\,\sC_{\Phi,0}(E). 
$$

More generally, for $\beta\in(-\Phi(-\be_n),\Phi(\be_n))$ Winterbottom in \cite{Winterbottom:1967} constructed an equilibrium shape of crystalls atop other material, which can be defined as 
\begin{equation}\label{truncated_ulffs}
W_\beta^\Phi:=\Omega\cap W^\Phi(-\beta\be_n),
\end{equation}
where $W^\Phi(z) = z + W^\Phi.$ As we have seen earlier, the ``half'' Wulff shape  $W_0^\Phi$ (which is also the Winterbottom shape with $\beta=0$) is not only an equilibrium, but also a global volume-constraint minimizer of $\sC_{\Phi,0}.$  The following result shows that this property is true also for other values of $\beta.$

\begin{theorem}\label{teo:main}
For any $\beta\in(-\Phi(-\be_n),\Phi(\be_n))$ 
\begin{equation}\label{min_capillary}
\inf_{E\in BV(\Omega;\{0,1\}),\,|E|= W_\beta^\Phi}\,\,\sC_{\Phi,\beta}(E) = \sC_{\Phi,\beta}(W_\beta^\Phi). 
\end{equation}
The equality holds if and only if $E = W^\Phi(b - \beta \be_n)$ for some $b\in\p\Omega.$ Equivalently, 
\begin{equation}\label{relat_isop_ineq}
\frac{\sC_{\Phi,\beta}(E)}{|E|^{\frac{n-1}{n} } } \ge \frac{\sC_{\Phi,\beta}(W_\beta^\Phi)}{|W_\beta^\Phi|^{\frac{n-1}{n} } },\qquad E\in BV(\Omega;\{0,1\}),
\end{equation}
and the equality holds if and only if $E = \Omega\cap(b - r\beta\be_n + rW^\Phi)$ for some $r>0$ and $b\in \p\Omega.$
\end{theorem}

Thus, the volume-constraint minimizers of $\sC_{\Phi,\beta}$ are precisely the horizontal translations of $W_\beta^\Phi.$ This result is well-known in the Euclidean case $\Phi=|\cdot|$ (see, for example, \cite[Theorem 19.21]{Maggi:2012}). To the best of my knowledge, there is no literature on the minimality of $W_\beta^\Phi$ except for cases where $\beta=0$ or $\Phi$ is Euclidean.

\subsection*{Acknowledgement} 
I thank Guido De Philippis and Francesco Maggi for useful discussions, especially, for showing their (unpublished) short notes on a generalization of Theorem \ref{teo:main} to more general cones with vertex at origin. I acknowledge support from the FWF Stand-Alone project P33716.

\section{Proof of Theorem \ref{teo:main}}

Owing \eqref{crs_2016} we provide an elementary proof of this result using only properties of anisotropies in $\R^n.$ We divide the prove into smaller steps. 
\smallskip

{\it Step 1: Introducing new anisotropies.} Fix any $\beta\in (-\Phi(-\be_n),\Phi(\be_n))$ and
$\eta^\pm\in\p\Phi(\pm \be_n),$ i.e.\footnote{Since a priori we are not assuming the regularity of $\Phi,$ there could be more than one possible choice of $\eta^\pm$. Moreover, since we are not assuming the evenness of $\Phi,$ in general we cannot claim $\eta^+=-\eta^-.$},
\begin{equation}\label{subdifferential}
\inner{\eta^\pm}{\pm\be_n} = \Phi(\pm \be_n)\quad \text{and}\quad \Phi^o(\eta^\pm) = 1,
\end{equation}
where $\p f$ is the subdifferential of a convex function $f.$ Consider the functions 
\begin{equation}\label{headlights}
\Psi_\beta(x):= 
\begin{cases}
\Phi(x) - \beta\,\inner{x}{\tfrac{\eta^+}{\Phi(\be_n)}} & \text{if $\beta\ge0,$}\\[2mm]
\Phi(x) + \beta\,\inner{x}{\tfrac{\eta^-}{\Phi(-\be_n)}} & \text{if $\beta<0,$}
\end{cases}
\qquad x\in \R^n,
\end{equation}
where for shortness we drop the dependence of $\Psi_\beta$ on $\Phi$ and the choice of $\eta^\pm.$ Notice  that such a technique of ``absorbing'' the relative adhesion coefficient into the anisotropy was already used in \cite{DPM:2015}. 

Let us show that $\Psi_\beta$ is an anisotropy in $\R^n$. Indeed, the convexity and positive one-homogeneity of $\Psi_\beta$ are obvious. 
Let us show that there exists $C_{\Psi_\beta} \ge c_{\Psi_\beta}>0$ such that 
\begin{equation}\label{norm_bounds_for_Psi}
c_{\Psi_\beta}|x|  \le \Psi_\beta(x) \le C_{\Psi_\beta}|x|,\quad x\in\R^n.
\end{equation}
Indeed, by \eqref{norm_bounds}
$$
\sup_{|x|=1} \Psi_\beta(x) \le \sup_{|x|=1} \Phi(x) + \tfrac{|\beta|\max\{|\eta^+|,|\eta^-|\}}{\Phi(\be_n)}\le C_\Phi +  \tfrac{|\beta|\max\{|\eta^+|,|\eta^-|\}}{\Phi(\be_n)} =:C_{\Psi_\beta}. 
$$
Thus, the second inequality in \eqref{norm_bounds_for_Psi} holds.

On the other hand, by the Young inequality\footnote{I.e., $\inner{x}{y}\le \Phi(x)\Phi^o(y)$ for all $x,y\in\R^n.$} and the second equality in \eqref{subdifferential} 
\begin{equation}\label{young_schal}
\inner{x}{\tfrac{\eta^\pm}{\Phi(\pm\be_n)}} \le \tfrac{\Phi(x)\Phi^o(\eta^\pm)}{\Phi(\pm\be_n)}  =  
\tfrac{1}{\Phi(\pm\be_n)}\, \Phi(x).
\end{equation}
Now if $\beta \ge0,$ then by  \eqref{young_schal} and \eqref{norm_bounds}
$$
\Psi_\beta(x) \ge \Big(1 -\tfrac{\beta}{\Phi(\be_n)}\Big)\,\Phi(x) \ge \tfrac{\Phi(\be_n) - \beta}{\Phi(\be_n)}\,c_\Phi|x|,
$$
and similarly, if $\beta<0,$
$$
\Psi_\beta(x) \ge \tfrac{\Phi(-\be_n) + \beta}{\Phi(-\be_n)}\,c_\Phi|x|.
$$
Thus, 
$$
c_{\Psi_\beta}:= c_\Phi \min\Big\{\tfrac{\Phi(\be_n) - |\beta|}{\Phi(\be_n)}, \tfrac{\Phi(-\be_n) - |\beta|}{\Phi(-\be_n)}\Big\}>0
$$
and the first inequality in \eqref{norm_bounds_for_Psi} holds. Therefore, $\Psi_\beta$ is an anisotropy in $\R^n.$
\smallskip

{\it Step 2: A representation of the capillary functional.} Let us show
$$
\sC_{\Phi,\beta} = P_{\Psi_\beta}(\cdot,\Omega).
$$
Indeed, since $\eta^\pm$ are constant, by the divergence theorem
$$
0 = \int_E \div\,\eta^\pm dx = \int_{\Omega\cap \p^*E} \inner{\eta^\pm}{\nu_E}\, d\cH^{n-1} - \int_{\p\Omega\cap \p^*E} \inner{\eta^\pm}{\be_n}\, d\cH^{n-1}. 
$$
Thus,
\begin{align*}
P_{\Psi_\beta}(E,\Omega) = & 
\int_{\Omega\cap\p^*E} \Phi(\nu_E)d\cH^{n-1} \mp \tfrac{\beta}{\Phi(\pm\be_n)} \int_{\Omega\cap \p^*E} \inner{\nu_E}{\eta^\pm}d\cH^{n-1}\\[2mm]
= & P_\Phi(E,\Omega) - \beta \,\tfrac{\inner{\eta^\pm}{\pm\be_n}}{\Phi(\pm\be_n)}\int_{\p\Omega}\chi_Ed\cH^{n-1}\\[2mm]
= & P_\Phi(E,\Omega) - \beta \int_{\p\Omega}\chi_Ed\cH^{n-1} = \sC_{\Phi,\beta}(E),
\end{align*}

{\it Step 3: Wulff shapes of $\Psi_\beta$ and $\Phi$.}  We claim
\begin{equation}\label{equal_wulffs}
W^{\Psi_\beta} = \mp\tfrac{\beta\eta^\pm}{\Phi(\pm \be_n)} + W^\Phi,
\end{equation}
where if $\beta\ge0$, we take "+" sign, otherwise we take "-" sign.  Indeed, assume that $\beta\ge0$ and take any $x$ with $\Psi_\beta^o(x) = 1,$ where $\Psi_\beta^o$ is the dual of $\Psi_\beta$ (since $\Psi_\beta$ is an anisotropy, its dual is well-defined and also is an anisotropy). We claim that 
\begin{equation}\label{carry_jimmo}
\Phi^o\Big(x + \tfrac{\beta\eta^+}{\Phi(\be_n)}\Big) = 1.
\end{equation}
Let $\xi \in \p\Psi_\beta^o(x),$ i.e., $\inner{x}{\xi} = 1$ and $\Psi_\beta(\xi) = 1.$ Then one can readily check that $x\in \p\Psi_\beta(\xi).$ Hence, using the explicit expression of $\Psi_\beta$ in \eqref{headlights} we can compute its subdifferential:
\begin{equation}\label{different_subdifferential}
\p\Psi_\beta(\theta) = \p\Phi(\theta) - \tfrac{\beta \eta^+}{\Phi(\be_n)} 
\end{equation}
at each $\theta\in\R^n\setminus \{0\},$ and get
$$
x = \zeta - \tfrac{\beta \eta^+}{\Phi(\be_n)}\quad\text{for some $\zeta\in \p\Phi(\xi)$}.
$$
Thus, 
$$
\Phi^o\Big(x + \tfrac{\beta \eta^+}{\Phi(\be_n)}\Big) = \Phi^o(\zeta) = 1. 
$$
On the other hand, if \eqref{carry_jimmo} holds, then $\zeta:=x + \tfrac{\beta \eta^+}{\Phi(\be_n)} \in \p \Phi(\xi)$ for some $\xi\ne0.$ This and \eqref{different_subdifferential} implies $x\in \p\Psi_\beta(\xi),$ i.e., $\Psi_\beta^o(x) = 1.$ Thus,
$$
\Phi^o\Big(x + \tfrac{\beta \eta^+}{\Phi(\be_n)}\Big) = 1\qquad\Longleftrightarrow\qquad \Psi_\beta^o(x)=1.
$$
Since both Wulff shapes are convex and their boundaries coincide, this implies \eqref{equal_wulffs}.

The case $\beta<0$ is analogous.
\smallskip

{\it Step 4: Translated Wulff shapes.}  Let us show that the translated $W^\Phi$ in \eqref{equal_wulffs} is a horizontal translation of truncated Wulf shapes $W_\beta^\Phi$ in \eqref{truncated_ulffs}. Indeed, consider the vector
$$
b: =  \mp\tfrac{\beta\eta^\pm}{\Phi(\pm \be_n)}  + \beta \be_n.
$$
By \eqref{subdifferential}
$$
\inner{b}{\be_n} = - \tfrac{\beta \inner{\eta^\pm}{\pm \be_n}}{\Phi(\pm\be_n)} + \beta = 0,
$$
and hence, $b\in\p\Omega.$ Therefore, the translated Wulff shape $\mp\tfrac{\beta\eta^\pm}{\Phi(\pm \be_n)} + W^\Phi$ is a horizontal translation of the translated Wulff shape $-\beta\be_n + W^\Phi.$
\smallskip

{\it Step 5: Minimality of truncated Wulff shape $W_\beta^\Phi$ in \eqref{truncated_ulffs}.} Applying \eqref{crs_2016} with $\Psi_\beta$ we find 
\begin{equation}\label{schedule_time}
\tfrac{P_{\Psi_\beta}(E,\Omega)}{|E|^{\frac{n-1}{n}}} \ge \tfrac{P_{\Psi_\beta}(W^{\Psi_\beta},\Omega)}{|\Omega\cap W^{\Psi_\beta}|^{\frac{n-1}{n}}},\qquad E\in BV(\Omega;\{0,1\}).
\end{equation}
The equality holds iff $E = b + rW^{\Psi_\beta}$ for some $r>0$ and $b\in\p\Omega.$ By steps 3 and 4, $W_\beta^\Phi$ is a horizontal translation of 
\begin{equation}\label{equalwulffs}
W^{\Psi_\beta} = W^\Phi(\mp \tfrac{\beta\eta^\pm}{\Phi(\pm \be_n)}) = b_0 + W^\Phi(-\beta\be_n)
\end{equation}
for some $b_0\in\p\Omega.$ In particular, we can use $W_\beta^\Phi$ in place of $\Omega\cap W^{\Psi_\beta}$ in \eqref{schedule_time}. Moreover, by step 2 $P_{\Psi_\beta}(\cdot,\Omega) = \sC_{\Phi,\beta},$ and hence, we can represent \eqref{schedule_time} as 
$$
\tfrac{\sC_{\Phi,\beta}(E)}{|E|^{\frac{n-1}{n}}} \ge \tfrac{\sC_{\Phi,\beta}(W_\beta^{\Phi})}{| W_\beta^{\Phi}|^{\frac{n-1}{n}}},\qquad E\in BV(\Omega;\{0,1\}),
$$
which is \eqref{relat_isop_ineq}. 
\smallskip

{\it Step 6: Conclusion of the proof of Theorem \ref{teo:main}.} 
Since $\p\Omega$ is a (horizontal) hyperplane, the set of all horizontal translations form an additive group. In particular, by step 5 and \eqref{equalwulffs} the sets
$$
E:=b + rW^{\Psi_\beta} = b + rW^\Phi(\mp \tfrac{\beta\eta^\pm}{\Phi(\pm \be_n)}) = (b + b_0r) + rW^\Phi(-\beta\be_n)
$$ 
are the only ones preserving the equality in \eqref{schedule_time}, or equivalently, the equality in \eqref{relat_isop_ineq} holds if and only if $E = \Omega\cap (b + W^\Phi(-\beta r\be_n))$ for some $b\in\p\Omega$ and $r>0.$ Finally, the assertions related to the equality \eqref{min_capillary} directly follows from \eqref{relat_isop_ineq}.

\end{document}